\documentclass[11pt, psamsfonts]{amsart}

\usepackage{amssymb,amsfonts}
\usepackage[all,arc]{xy}
\usepackage{enumerate}
\usepackage{mathrsfs}
\usepackage{thmtools}
\usepackage{thm-restate}
\usepackage{datetime}
\usepackage{hyperref}
\usepackage{graphicx}
\usepackage{cleveref}
\usepackage{comment}
\usepackage{dirtytalk}
\usepackage{xcolor}
\definecolor{amethyst}{rgb}{0.6, 0.4, 0.8}
\definecolor{darkspringgreen}{rgb}{0.09, 0.45, 0.27}
\newenvironment{bnote}{\color{blue}}{}

\newenvironment{gnote}{\color{darkspringgreen}}{}
\newenvironment{rnote}{\color{red}}{}

\calclayout

\newtheorem{thm}{Theorem}[section]
\newtheorem{cor}[thm]{Corollary}
\newtheorem{prop}[thm]{Proposition}
\newtheorem{lem}[thm]{Lemma}
\newtheorem{conj}[thm]{Conjecture}
\newtheorem{ques}[thm]{Question}

\theoremstyle{definition}
\newtheorem{defn}[thm]{Definition}

\theoremstyle{remark}
\newtheorem{rem}[thm]{Remark}

\DeclareMathOperator{\cart}{\square}

\makeatletter
\let\c@equation\c@thm
\makeatother
\numberwithin{equation}{section}

\makeatletter
\newcommand*\bigcdot{\mathpalette\bigcdot@{.5}}
\newcommand*\bigcdot@[2]{\mathbin{\vcenter{\hbox{\scalebox{#2}{$\m@th#1\bullet$}}}}}
\makeatother

\makeatletter
\def\subsection{\@startsection{subsection}{3}%
  \z@{.5\linespacing\@plus.7\linespacing}{.1\linespacing}%
  {\bfseries}}
\makeatother

\newcommand{\N} { \mathbb{N}}

\newcommand{\Z} { \mathbb{Z}}

\newcommand{\cut} { \backslash}

\usepackage{etoolbox}
\makeatletter
\patchcmd{\@maketitle}
  {\ifx\@empty\@dedicatory}
  {\ifx\@empty\@date \else {\vskip3ex \centering\footnotesize\@date\par\vskip1ex}\fi
   \ifx\@empty\@dedicatory}
  {}{}
\patchcmd{\@adminfootnotes}
  {\ifx\@empty\@date\else \@footnotetext{\@setdate}\fi}
  {}{}{}
\makeatother

\title[Unimodality and Monotonic Portions of Certain Domination Polynomials]{Unimodality and Monotonic Portions of Certain Domination Polynomials}
\author{Amanda Burcroff}
\address{Department of Mathematical Sciences \\
Durham University\\
Durham, United Kingdom DH1 3LE}
\email{amanda.g.burcroff@durham.ac.uk \textrm{ or } aburcroff@math.harvard.edu}
\author{Grace O'Brien}
\address{Department of Mathematics\\ 
University of Michigan\\
Ann Arbor, Michigan 48009}
\email{graceob@umich.edu}

\begin{document}
\maketitle

\begin{abstract}
    Given a simple graph $G$ on $n$ vertices, a subset of vertices $U \subseteq V(G)$ is \emph{dominating} if every vertex of $V(G)$ is either in $U$ or adjacent to a vertex of $U$.  The \emph{domination polynomial} of $G$ is the generating function whose coefficients are the number of dominating sets of a given size.  We show that the domination polynomial is unimodal, i.e., the coefficients are non-decreasing and then non-increasing, for several well-known families of graphs.  In particular, we prove unimodality for spider graphs with at most $400$ legs (of arbitrary length), lollipop graphs, arbitrary direct products of complete graphs, and Cartesian products of two complete graphs. We show that for every graph, a portion of the coefficients are non-increasing, where the size of the portion depends on the upper domination number, and in certain cases this is sufficient to prove unimodality.  Furthermore, we study graphs with $m$ \emph{universal vertices}, i.e., vertices adjacent to every other vertex, and show that the
    last $(\frac{1}{2} - \frac{1}{2^{m+1}}) n$ coefficients of their domination polynomial are non-increasing.  
\end{abstract}

\section{Introduction}
Graph domination has become a mainstream branch of graph theory over the last half century, but only in the past decade have researchers asked enumerative questions about the number of dominating sets of a particular size in a fixed graph.  In this paper, we study these quantities on several families of graphs and determine special properties of these counts as sequences, iterating over the size of the dominating sets for a given graph.

Let $G = (V,E)$ be a simple graph on $n$ vertices.  A set $U \subseteq V$ is called \emph{dominating} if every vertex of $G$ is either in $U$ or adjacent to a vertex of $U$.  Let $d_i(G)$ denote the number of dominating sets of a graph $G$ of a fixed size $i$. A significant portion of the study of domination in graphs has focused on finding the minimum size of a dominating set in $G$, the \textit{domination number}, $\gamma(G)$. Note that the minimum $i$ such that $d_i(G) > 0$ is precisely $\gamma(G)$, so calculating these quantities $d_i(G)$ encompasses the task of determining the domination number. It is often convenient to look at these $d_i(G)$ in terms of generating functions; the {\it domination polynomial}, $D(G,x)$, is given by 
$$D(G,x) = \sum_{i=1}^n d_i(G)x^i\,.$$ 

Even without precisely calculating the quantities $d_i(G)$, one can investigate properties of their distribution, such as their relative order. Alikhani and Peng posed the following conjecture in 2014:
\begin{conj}[\cite{AP}]\label{conj:unimodal}
The domination polynomial of any graph is unimodal.
\end{conj}
Note that, in general, the coefficients of the domination polynomial are not log-concave (see an example of minimum size on $9$ vertices in \cite{BeBr}).  This implies that the domination polynomials are not real-rooted, and in fact it is known that the roots are dense in the complex plane \cite{BrTu}. 

There has been an assortment of partial progress toward Conjecture \ref{conj:unimodal} (see \cite{Bea} for a survey of recent results). In 2020, Beaton and Brown \cite{BeBr} show that paths, cycles, and complete multipartite graphs all have unimodal domination polynomials. Moreover, they prove that graphs with sufficiently large minimum degree have unimodal domination polynomials (see Theorem \ref{thm:mindeg}), which they use to show that almost all graphs have unimodal domination polynomials.  In \cite{AJ}, Alikhani and Jahari demonstrate unimodality of the domination polynomials for several families of graphs, including every friendship graph as well as the corona of any graph with $P_3$ or $K_n$.  

In this paper, we identify other families of graphs with unimodal domination polynomials and study the behavior of certain portions of the coefficients. The families of graphs that we study generally have certain extremal degree properties, i.e., have many vertices of low degrees, high degrees, or equal degrees.  We start with spider graphs, which are trees with only one vertex of degree larger than $2$.  We extend our methods used for spider graphs to handle lollipop graphs, which consist of vertices of degree at most $2$ and all other vertices of high degree.  We then move on to two well-structured families of regular graphs, namely the direct and Cartesian products of complete graphs.  We furthermore show that graphs with a very low upper domination number have unimodal domination polynomials; such graphs generally have all vertices of large degree.  Lastly, we focus on graphs with one or more universal vertices, which have the maximum possible degree.  Loosely speaking, this paper is structured so that the degree sequences of the graphs involved are lexicographically increasing (though, of course, such an ordering depends on the choice of several parameters).

More explicitly, the remainder of this paper is structured as follows.  We begin by providing the necessary background in Section \ref{sec: prelim}.  In Section \ref{sec:spiders}, we consider spider graphs having a bounded number of legs, each of arbitrary length.  In particular, we prove that all spider graphs with at most $400$ legs have unimodal domination polynomials.  We remark that the bound $400$ is somewhat arbitrary; a portion of the proof relies on a number of computer checks that could be carried out further, but in our case slowed down around this threshold. Our work includes a generalization of the methods developed by Beaton and Brown \cite{BeBr} for paths.  We furthermore show that any lollipop graph, i.e., a graph obtained by appending a path to a complete graph, has a unimodal domination polynomial in Section \ref{sec:lollipop}.  

Next, we investigate connected direct and Cartesian products of complete graphs. The interplay between dominating sets and graph products has been the subject of myriad studies, as surveyed in \cite[Chapter 28]{HIK} and \cite{NR}.  In general, it can be hard to estimate how the sequence $d_i(G)$ behaves under graph products.  One of the most famous open problems in graph domination is Vizing's conjecture that the domination number of a Cartesian product of graphs is at least the product of their domination numbers \cite{Viz}.  In Sections \ref{sec:direct} and \ref{sec:cartesian}, we prove the unimodality of the direct product of arbitrarily many complete graphs and the Cartesian product of two complete graphs, respectively. The former family includes the unitary Cayley graphs of $\Z/q\Z$ for every squarefree integer $q$, which have received recent attention for their extremal domination properties \cite{Bur, Bur2, DI, Vem}.

In Section \ref{sec:upper dom}, we use a Hall-type argument to show that a portion of the domination polynomial coefficients are non-increasing, where the size of this portion depends on the upper domination number.  The \emph{upper domination number} $\Gamma(G)$ of a graph $G$ is the maximum size of a minimal dominating set in $G$.  Our bound in particular implies that graphs of odd order with $\Gamma(G) \leq 4$ and graphs of even order with $\Gamma(G) \leq 3$ have unimodal domination polynomials.

In Section \ref{sec:universal}, we study graphs with $m$ universal vertices, as suggested in the open problems section of \cite{BeBr}.  A \emph{universal vertex} is a vertex that is adjacent to every other vertex, i.e., has degree $n-1$.  We show that the last $\left(\frac{1}{2}-\frac{1}{2^{m+1}}\right)n$ coefficients of the domination polynomial are non-increasing.  As a direct corollary, any graph with at least $\log_2(n)-1$ universal vertices has a unimodal domination polynomial.  We conclude with several directions for future research in Section \ref{sec: further directions}.





\section{Preliminaries}\label{sec: prelim}
Here we provide some basic definitions, an overview of the relevant families of graphs and graph operations, as well as a collection of previous results.

A polynomial is called \emph{unimodal} if its coefficients are non-decreasing and then non-increasing. We order the coefficients of a polynomial $f(x)$ of degree $n$ by increasing order of the corresponding powers of $x$; thus when we refer to the last $k$ coefficients, we mean the coefficients of $x^{n-k+1}$, $x^{n-k+2}$, \ldots, $x^n$.  If $f(x)$ is unimodal, a coefficient of $x^k$ is a \emph{mode} if it achieves the maximum value over all coefficients; in this case, we say $f$ \emph{has a mode at} $k$.  Note that a unimodal polynomial may have multiple modes.  Given two unimodal polynomials $f(x)$ and $g(x)$, we refer to the \emph{distance} between modes appearing as the coefficients of $x^{k_f}$ and $x^{k_g}$ in $f(x)$ and $g(x)$, respectively, as the quantity $|k_f - k_g|$.

We begin by discussing the domination polynomials of spiders.  The \emph{spider graph} (or \emph{spider}) $S(\lambda_1,\ldots,\lambda_t)$ is the graph formed by taking $t$ paths, each on $\lambda_i$ vertices, and connecting one end of each path to an additional central vertex.  The induced path graph on a set of vertices comprising one of the original paths is called a \emph{leg of length $\lambda_i$}.

\begin{figure}[ht!]
    \centering
    \includegraphics[width=.3\linewidth]{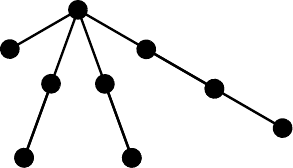}
    \caption{The spider graph $S(1,2,2,3).$}
    \label{fig:spider}
\end{figure}

In the study of spiders, we make use of the several common graph notions, defined here. Let $G = (V,E)$ be a graph and $v$ be a vertex.  The \emph{open neighborhood of $v$}, denoted by $N(v)$, is the set of all neighbors of $v$.  The \textit{closed neighborhood of $v$}, denoted by $N[v]$, is the open neighborhood of $v$ along with $v$ itself, i.e., $N[v] = \{u \in V : \{u,v\} \in E\} \cup \{v\}$.   We also look at subgraphs obtained by the deletion of one or more vertices.  For any $U \subseteq V$, the graph $G - U$ is the graph obtained by removing all vertices in $U$ and any edges incident to a vertex in $U$.  When $U = \{v\}$ consists of a single vertex, we often write $G - v$ for $G - \{v\}$.  The subgraph $G/v$ obtained by \emph{contracting} the vertex $v$ has vertex set $V\cut \{v\}$, and its edge set consists of all edges in $G-\{v\}$ plus all edges between two vertices in $N(v)$.

After our discussion of spiders, we move on to graph products. There are multiple ways to define a product of graphs; we focus on two of the most common, namely the direct and Cartesian products. Let $G = (V(G), E(G))$ and $H = (V(H),E(H))$ be graphs. The \textit{direct product} (or \emph{tensor product}) of $G$ and $H$ is the graph $G \times H$ on $V(G\times H) = V(G) \times V(H)$ where $(g,h)$ and $(g',h')$ are adjacent if $\{g,g'\} \in E(G)$ and $\{h,h'\} \in E(H)$. The \textit{Cartesian product} (or \emph{box product}) is the graph $G \cart H$ on the vertices of $V(G) \times V(H)$ where $(g,h)$ and $(g',h')$ are adjacent if $g=g'$ and $\{h,h'\} \in E(H)$, or $\{g,g'\} \in E(G)$ and $h = h'$.  When we take the product of multiple graphs simultaneously, the definitions follow associatively from the definitions for the product of two graphs.


Finally, we introduce two theorems, which we utilize several times throughout this paper.
\begin{thm}[\protect{\cite[Theorem 3.2]{BeBr}}]\label{thm:mindeg}
If $G$ is a graph on $n$ vertices with minimum degree 
\\ $\delta(G) \geq 2\log_2(n)$, then $D(G,x)$ is unimodal with a mode at $\left \lceil{\frac{n}{2}}\right \rceil$.
\end{thm}

\begin{thm}[\protect{\cite[Theorem 5]{AP}}]\label{thm: first half}
Let $G$ be a graph on $n$ vertices.  Then for every $0 \leq i < \frac{n}{2}$, we have $d_i(G) \leq d_{i+1}(G)$.
\end{thm}

Theorem \ref{thm:mindeg} is the foundation for Beaton and Brown's proof that almost all graphs (in the context of the Erd{\H o}s-R{\'e}nyi graph model) are unimodal (\cite[Theorem 3.3]{BeBr}).  Their result allows us to simplify the consideration of many families of graphs. The proof of Theorem \ref{thm: first half} is a direct application of Hall's Theorem (see Theorem \ref{thm: Hall}). In general, making statements about the first half of the coefficients of a domination polynomial has been relatively straightforward, while claims about the second half have proved more difficult.

\section{\texorpdfstring{Spiders With at Most $400$ Legs}{Spiders With at Most 400 Legs}}\label{sec:spiders}

Trees are a particularly elusive case of Conjecture \ref{conj:unimodal} due to their low minimum degree, as noted in \cite{BeBr}. As a first step towards understanding the behavior of the domination polynomials for trees in general, we study spider graphs, a family of trees. In this section, we show that all spider graphs with at most $400$ legs have unimodal domination polynomials.

In order to establish unimodality for the domination polynomials of spiders, we adapt methods developed by Beaton and Brown to show that $D(P_n,x)$ is unimodal for all $n$, where $P_n$ is the path on $n$ vertices.  In particular, we require the following theorem, which allows us to establish unimodality for a sequence of polynomials.

\begin{thm}\label{thm: BeBr poly recursion}\emph{(\cite[Theorem 2.2]{BeBr})}
Suppose  we  have  a  sequence  of  polynomials $(f_n)_{n \geq 1}$ with  non-negative  coefficients that satisfy $f_n=x(f_{n-1}+f_{n-2}+f_{n-3})$
for $n\geq 4$.  Let $\mathcal{P}_n$ denote the following property: for all $i \in \{1, \ldots , n\}$, $f_i$ is unimodal with a mode at $\mu_i$ and for $i\geq 2$, we additionally have $0\leq \mu_i-\mu_{i-1}\leq 1$.  Assume $\mathcal{P}_4$ holds.  Then $\mathcal{P}_n$ holds for all $n\geq 1$ (and so each $f_n$ is unimodal).
\end{thm}

In order to apply Theorem \ref{thm: BeBr poly recursion} to spiders, we use the following recurrence due to Kotek, Preen, Simon, Tittmann, and Trinks.

\begin{prop}\emph{(\cite[Proposition 3.1]{KPS})} \label{prop: recursion}
Let $G$ be a graph. If there exist $u, v \in V(G)$ such that $N[u] \subseteq N[v]$, then $$D(G,x) = xD(G/v) + D(G-v) + xD(G-N[v])\,.$$
\end{prop}

We show that if $G$ can be formed by appending a path on $3$ vertices to a graph, then we can confirm the unimodality of $D(G,x)$ by checking certain unimodality and mode conditions on subgraphs of $G$.  We now proceed to characterize which subgraphs must be considered.

Fix a graph $H$ on $n$ vertices and a vertex $v \in V(H)$.  Let $H_v^{(\ell)}$ denote the graph obtained by appending a path on $\ell$ vertices by adding an edge between $v$ and a leaf, i.e., a vertex of degree $1$, of the path. We set $H_v^{(0)} := H$.  Thus, $H_v^{(\ell)}$ has $n + \ell$ vertices for all $\ell \geq 0$.  

\begin{figure}[ht!]
    \centering
    \includegraphics[width=.9\linewidth]{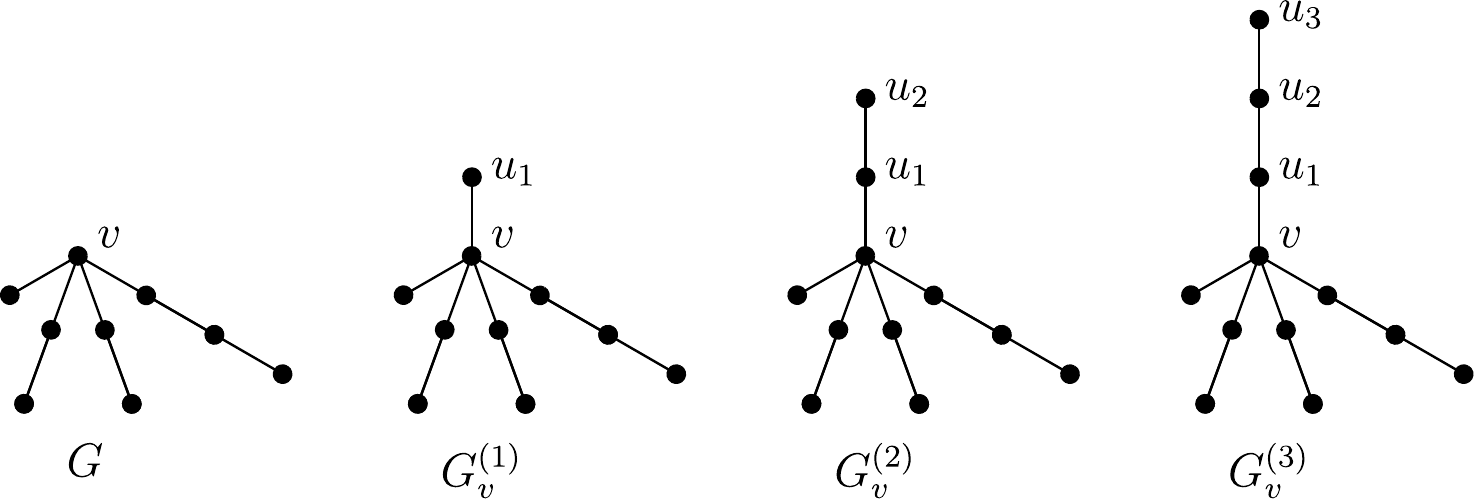}
    \caption{Illustrations of the graphs constructed in Lemma \ref{lem: path recursion}, where 
    \\$G = S(1,2,2,3)$ and $v$ chosen to be the unique vertex of degree at least $3$.}
    \label{fig: path recursion example}
\end{figure}

\begin{lem}\label{lem: path recursion}
Let $G = (V,E)$ be a graph, and fix $v \in V(G)$. Suppose we have that
\begin{enumerate}[(i)]
    \item $D(G,x)$, $D(G_v^{(1)},x)$, $D(G_v^{(2)},x)$, and $D(G_v^{(3)},x)$ are unimodal with modes at $\mu_0$, $\mu_1$, $\mu_2$, and $\mu_3$, respectively, and
    \item $0 \leq \mu_i - \mu_{i-1} \leq 1$ for $i \in \{1,2,3\}$.
\end{enumerate}
Then $D(G_v^{(\ell)},x)$ is unimodal for all $\ell \in \N$.
\end{lem}
\begin{proof}
Let $u_i$ be the unique vertex in $V\left(G_v^{(i)}\right) \cut V\left(G_v^{(i-1)}\right)$ for $i \in \{1,2,3\}$, as illustrated in Figure \ref{fig: path recursion example}.  Observe that $G_v^{(3)}/u_2 \cong G_v^{(2)}$ and $G_v^{(3)} - N[u_2] = G_v^{(3)} - \{u_1,u_2,u_3\} = G$. 

Since $N[u_3] \subseteq N[u_2]$ in $G_v^{(3)}$, by Proposition \ref{prop: recursion} we have
$$D\left(G_v^{(3)},x\right) = xD\left(G_v^{(2)},x\right) + D\left(G_v^{(3)}- \{u_2\},x\right) + xD\left(G,x\right).$$
Moreover, note that  $G_v^{(3)} - \{u_2\} = G_v^{(1)} \cup \{u_3\}$.  Since $u_3$ must be included in every dominating set of $G_v^{(1)} \cup \{u_3\}$, we have $D\left(G_v^{(3)} - \{u_2\},x\right) = xD\left(G_v^{(1)},x\right)$.  Thus, Theorem \ref{thm: BeBr poly recursion} applies with $f_n = D\left(G_v^{(3)},x\right)$, $f_{n-1} = D\left(G_v^{(2)},x\right)$, $f_{n-2} = D\left(G_v^{(1)},x\right)$, and $f_{n-3} = D\left(G,x\right)$.
\end{proof}

While we use Lemma \ref{lem: path recursion} to show that spiders with a bounded number of legs have unimodal domination polynomials, similar techniques could be applied to other families of graphs satisfying similar hypotheses.  We need only check that these conditions hold when attaching paths on at most $3$ vertices to establish unimodality after attaching paths of arbitrary lengths.  Once such family where similar methods can be used is the set of lollipop graphs, which we investigate in Section \ref{sec:lollipop}.  

We now establish the precise domination polynomial of a spider graph with legs of length at most $3$.  This classification is then used to complete the finite number of checks needed to ensure the conditions of Lemma \ref{lem: path recursion} hold for spiders with a bounded number of legs, each with length at most $3$.

\begin{lem}\label{lem: spider poly}
Let $S$ be a spider comprised of $\ell_1$ legs of length $1$, $\ell_2$ legs of length $2$, and $\ell_3$ legs of length $3$.  If $\ell_1 > 0$, we have
$$D(S,x) = x(1+x)^{\ell_1}(2x+x^2)^{\ell_2}(2x+3x^2 + x^3)^{\ell_3} + x^{\ell_1}(2x + x^2)^{\ell_2}(x+3x^2+ x^3)^{\ell_3}\,.$$
If $\ell_1 = 0$, we have
$$D(S,x) = x(2x + x^2)^{\ell_2}(2x +3x^2 + x^3)^{\ell_3} + (2x + x^2)^{\ell_2}(x + 3x^2 + x^3)^{\ell_3} - x^{\ell_2}(x + x^2)^{\ell_3}\,.$$
\end{lem}
\begin{proof}
Let $z$ denote the center vertex of the spider.  We claim that 
$$A(x) := x(1 + x)^{\ell_1}(2x + x^2)^{\ell_2}(2x + 3x^2 + x^3)^{\ell_3}$$
is the generating function whose coefficients count the number of dominating sets of $S$ that include $z$.  Observe that $2x + 3x^2 + x^3$ (resp. $2x+x^2$, $1 + x$) counts the number of sets of a path of length $3$ (resp. $2$, $1$) whose closed neighborhood includes all vertices except possibly one endpoint of the path.  The dominating sets of $S$ are formed by taking such a subset of vertices in each of the $\ell_i$ paths of length $i$, and then attaching the (possibly uncovered) endpoint of each path to $z$ via an edge.   

We now proceed to count the number of dominating sets of $S$ that do not contain $z$.  Observe that $x + 3x^2 + x^3$  (resp. $2x + x^2$, $x$) counts the number of dominating sets of a path of length $3$ (resp. $2$, $1$).  Suppose that we choose a subset of vertices $W \subset V(S)$ formed by taking such a subset of vertices in the $\ell_i$ paths of length $i$, and then attaching an endpoint of each path to $z$ via an edge.  Note that there are $x^{\ell_1}(2x + x^2)^{\ell_2}(x + 3x^2 + x^3)^{\ell_3}$ such sets.  Our construction of $W$ guarantees that $W$ dominates all vertices of $S$ except for possibly $z$; we now consider two cases regarding when $z$ is dominated.
\begin{enumerate}[(i)]
\item If $\ell_1 > 0$, then one of the legs of length $1$ in $S$ must have its leaf included in $W$.  Thus, $W$ dominates $z$. So, by the above, $W$ is a dominating set of $S$.
\item If $\ell_1 = 0$, then $W$ is non-dominating precisely when all the endpoints of the paths that are attached to $z$ are not included in $W$.  Such sets are counted by $x^{\ell_2}(x + x^2)^{\ell_3}$.
\end{enumerate}
We can therefore see that the number of dominating sets of $S$ not containing $z$ is given by the coefficients of
$$B(x) :=
\begin{cases}
x^{\ell_1}(2x + x^2)^{\ell_2}(x + 3x^2 + x^3)^{\ell_3} & \text{ if } \ell_1 > 0\,,\\
(2x + x^2)^{\ell_2}(x + 3x^2 + x^3)^{\ell_3} - x^{\ell_2}(x + x^2)^{\ell_3} & \text{ if } \ell_1 = 0\,.
\end{cases}
$$
Setting $D(S,x) = A(x) + B(x)$ yields the desired equalities.
\end{proof}

\begin{lem}\label{lem: spider recursion}
Let $S = S(\lambda_1,\ldots,\lambda_t)$ be a spider with $t \leq 400$, and fix $i$ such that $\lambda_i > 3$.  Let $S' = S(\lambda_1,\ldots,\lambda_{i-1},\lambda_i+1,\lambda_{i+1},\ldots,\lambda_t)$.  If $D(S,x)$ is unimodal with a mode at $\mu$, then $D(S',x)$ is unimodal with a mode either at $\mu$ or $\mu+1$. 
\end{lem}
\begin{proof}
Apply Lemma \ref{lem: path recursion} to the vertices at the end of the leg of length $\lambda_i$; in particular, set $u_3$ to be the leaf of this leg, $u_2$ the unique neighbor of $u_3$, $u_1$ the other neighbor of $u_2$, and $v$ the other neighbor of $u_1$.  
\end{proof}

\begin{thm}
All spider graphs with at most $400$ legs (each of arbitrary length) have unimodal domination polynomials.
\end{thm}
\begin{proof}
We first note that for $t \leq 400$, one can check via computer (as we have successfully done) that a spider $S$ with $t$ legs, each of length at most $3$, satisfies the hypotheses of Lemma \ref{lem: path recursion}.  The domination polynomials for these spiders are given explicitly by Lemma \ref{lem: spider poly}, so the computer check is straightforward.  It then follows immediately from Lemma \ref{lem: spider recursion} that every spider with at most $400$ legs has a unimodal domination polynomial.
\end{proof}

The natural next step is to determine if $D(S,x)$ is unimodal for any spider $S$. Note that the unimodality of these polynomials is implied by Alikhani and Peng's conjecture (\ref{conj:unimodal}). However, proving the following conjecture on spiders with legs of length at most $3$ would also be sufficient to prove unimodality for all spiders by applying Lemma \ref{lem: path recursion}. 

\begin{conj}\label{conj: spider modes}
Let $S$ be a spider with legs of length at most $3$, and $S'$ be a spider obtained by deleting a leaf from $S$.  Then $D(S,x)$ and $D(S',x)$ are unimodal with at least one pair of modes at $\mu_S$ and $\mu_{S'}$, respectively, such that $|\mu_S - \mu_{S'}| \leq 1$.
\end{conj}

\section{Lollipop Graphs}\label{sec:lollipop}
Furthering our study of graphs with low minimum degree, we show that all lollipop graphs have unimodal domination polynomials.  We use recursive methods with respect to appending paths, similar to those used for spider graphs in the previous section.  

\begin{defn}
The \emph{$(m,n)$-lollipop graph} $L_{m,n}$ is the graph on $m + n$ vertices consisting of a complete graph on $m$ vertices and a path graph on $n$ vertices, connected at a leaf of the path via an edge.
\end{defn}

\begin{lem}\label{lem: lollipop cases}
For $1 \leq n \leq 3$ and $m\geq 3$, we have that $D(L_{m,n},x)$ is unimodal with a mode at $\mu$, where
\[
\mu =
\begin{cases}
\left\lfloor\frac{m}{2}\right\rfloor + 1 & \text{ if } n = 1\,,\\
\left\lceil\frac{m}{2}\right\rceil + 1 & \text{ if } n = 2\,,\\
\left\lfloor\frac{m}{2}\right\rfloor + 2 & \text{ if } n = 3\,.\\
\end{cases}
\]
\end{lem}
\begin{proof}
We break into cases depending on the value of $n$.  Beaton and Brown \cite{BeBr} completed the case $m=3$ in the remark after their Corollary 2.3, so we can take $m \geq 4$.  Thus, $m + n \geq 7$, so in particular, we know by Theorem \ref{thm: first half} that the first $4$ coefficients are non-decreasing.\\\\
{\bf Case $n = 1$:} Note that a set of at least $2$ vertices in $L_{m,1}$ is non-dominating if and only if it does not contain the leaf or the leaf's unique neighbor, i.e., the set is contained in the remaining $m-1$ vertices.  Hence, for $i \geq 2$ we have
$$d_i(L_{m,1}) = {m + 1 \choose i} - {m - 1 \choose i} = {m \choose i - 1} + {m - 1 \choose i - 1}\,.$$
The first summand ${m \choose i -1}$ has modes at $\left\lceil\frac{m}{2}\right\rceil + 1$ and $\left\lfloor\frac{m}{2}\right\rfloor + 1$, while the second summand ${m-1 \choose i-1}$ has modes at $\left\lceil\frac{m-1}{2}\right\rceil + 1$ and $\left\lfloor\frac{m-1}{2}\right\rfloor + 1$.  In particular, both sequences have a mode at 
$\left\lfloor\frac{m}{2}\right\rfloor + 1 = \left\lceil\frac{m - 1}{2}\right\rceil  + 1$.
Hence, $D(L_{m,1},x)$ is unimodal with a mode at the desired position.\\\\
{\bf Case $n = 2$:} A set of at least $3$ vertices in $L_{m,2}$ is dominating if and only if it nontrivially intersects the closed neighborhood of the unique leaf.  Thus, for $i \geq 3$ we have
$$d_i(L_{m,2}) = {m + 2 \choose i} - {m \choose i} = {m + 1 \choose i -1} + {m + 1 \choose i} - {m \choose i} = {m + 1 \choose i -1} + {m \choose i-1}\,.$$
Similar to our analysis in Case 1, though with the indexing for $m$ shifted by one, we can see that $D(L_{m,2},x)$ is unimodal with a mode at $\left\lfloor\frac{m + 1}{2}\right\rfloor + 1 = \left\lceil\frac{m}{2}\right\rceil  + 1$.

\begin{figure}[ht!]
    \centering
    \includegraphics[width=.25\linewidth]{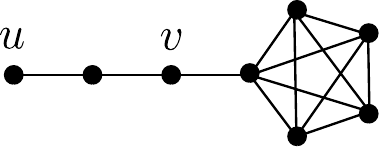}
    \caption{The lollipop $L_{5,3}$ with designated vertices $u$ and $v$.}
    \label{fig:lollipop}
\end{figure}
\noindent {\bf Case $n = 3$:} Let $u$ denote the unique leaf of $L_{m,3}$, and let $v$ be the unique vertex of distance $2$ from $u$ (see Figure \ref{fig:lollipop}). A set of at least $4$ vertices is dominating if and only if it dominates both $u$ and $v$.  This is because the unique neighbor of $u$ is dominated provided that $u$ is dominated, and moreover, by the pigeonhole principle at least one vertex of $K_m$ is included in the set hence all vertices of $K_m$ are dominated.  There are ${m + 1 \choose i}$ sets of size $i$ that do not dominate $u$, ${m \choose  i}$ sets of size $i$ that do not dominate $v$, and ${m - 1 \choose i}$ sets of size $i$ that dominate neither $u$ nor $v$.  Hence, we have
\begin{align*}
    d_i(L_{m,3}) &= {m + 3 \choose i}  - {m + 1 \choose i} - {m \choose i} + {m - 1 \choose i}\\
    &= {m + 2 \choose i - 1} + {m + 2 \choose i} - {m + 1 \choose i} - {m - 1 \choose i - 1}\\
    &= {m + 2 \choose i - 1} + {m \choose i-2} + {m \choose i-1} - {m-1 \choose i-1}\\
    &= {m + 2 \choose i -1} + {m \choose i-2} + {m-1 \choose i-2}\,.
\end{align*}

Similar to the previous two cases, we have that ${m \choose i-2} + {m-1 \choose i-2}$ is unimodal with a mode at $\left\lfloor \frac{m}{2}\right\rfloor + 2$.  Moreover, ${m + 2 \choose i-1}$ has a mode at $\left\lfloor \frac{m + 2}{2}\right\rfloor + 1 = \left\lfloor \frac{m}{2}\right\rfloor + 2$.  Hence, $D(L_{m,3},x)$ is as desired.
\end{proof}

\begin{thm}
All lollipop graphs have unimodal domination polynomials.
\end{thm}
\begin{proof}
This follows directly from Lemma \ref{lem: path recursion} and Lemma \ref{lem: lollipop cases}.
\end{proof}

\section{Products of Complete Graphs}\label{sec: prod}
While the minimum size of a dominating set in a product of graphs has been studied for decades \cite{HIK, NR, Viz}, there has been little work toward understanding the dominating sets of larger sizes in a graph product.  In general, this relationship can be rather complex, so we focus on a family of graphs that has been shown to have interesting domination properties in other contexts: direct and Cartesian products of complete graphs.  The direct
 product of complete graphs has been shown to have extremal properties in the domination chain \cite{Bur,DI,Vem}, while the Cartesian products of two complete graphs were used by the first author to demonstrate the tightness of Bre{\v s}ar, Klav{\v z}ar, and Rall's inequality $\Gamma(G\times H) \geq \Gamma(G)\Gamma(H)$, which holds for any graphs $G$ and $H$ \cite{Bur2, BKR}.

\subsection{Direct Products of Complete Graphs} \label{sec:direct}

In this section, we investigate the domination polynomials of the connected direct product of complete graphs. We first prove unimodality for certain regular graphs in Lemma \ref{lem:regular}, including graphs of the form $K_2 \times K_n$. Then, we consider the product of an arbitrary number of complete graphs.

\begin{lem} \label{lem:regular}
Let $G$ be an $m$-regular graph on $2n$ vertices for some $3 \leq n-1 \leq m < 2n$. Then, $D(G,x)$ is unimodal with a mode at $n$.
\end{lem}




\begin{proof}
Let $G=(V,E)$ be a graph as described. First, we notice that by Theorem \ref{thm: first half}, we have $d_1(G) \leq d_2(G) \leq \cdots \leq d_{n-1}(G) \leq d_n(G)$.  By the degree condition, each vertex $v$ has a closed neighborhood of size at least $n$. Every set of vertices that does not dominate $v$ must be a subset of $V(G)\cut N[v]$, which has size $2n - |N[v]| \leq n$.  Iterating over the vertices, we see that there are at most $2n$ non-dominating sets of size $n$, and furthermore every set of size at least $n+1$ is dominating.  Thus, ${2n \choose n} - 2n \leq d_n(G) \leq {2n \choose n}$, and $d_{n+r}(G) = {2n \choose n+r}$ for $1 \leq r \leq n$.  In particular, we clearly have $d_{n+1}(G) \geq d_{n+2}(G) \geq \cdots \geq d_{2n}(G)$.  For $n \geq 4$, it is straightforward to check that 
$${2n \choose n} - 2n \geq {2n \choose n+1}\,.$$
Since $d_{n-1}(G) \leq {2n \choose n-1} = {2n \choose n+1} = d_{n+1}(G) \leq d_n(G)$, we can conclude that $D(G,x)$ is unimodal with a mode at $n$.  
\end{proof}

\begin{rem}
Lemma \ref{lem:regular} can also be proven using the minimum degree argument from Theorem \ref{thm:mindeg} and a finite number of exceptional cases.  However, the exceptional cases would include checking the domination polynomials of all $(n-1)$-regular graphs on $2n$ vertices for $n \leq 9$, which appears computationally challenging. 
\end{rem}

We can now focus on connected direct products of at least three complete graphs. Note that a direct product of complete graphs is connected if and only if each complete graph has order at least $2$ and at most one of the complete graphs has order exactly $2$.  We begin by proving a technical lemma, which we will use later to show that all but finitely many direct products of complete graphs satisfy the hypotheses of Theorem \ref{thm:mindeg}.

\begin{lem}\label{lem: cart product ineq}
Suppose we have positive integers $2 \leq n_1 \leq n_2 \leq \cdots \leq n_t$ with $n_2 \geq 3$; and additionally if $t = 2$, then $n_1 \geq 3$.  If there exist positive integers $\ell_i \leq n_i$ for $1 \leq i \leq t$ such that
$$(\ell_1 - 1)(\ell_2 - 1)\cdots (\ell_t - 1) \geq 2\log_2(\ell_1\ell_2\cdots\ell_t)\,,$$
then 
$$(n_1 - 1)(n_2 - 1)\cdots (n_t - 1) \geq 2\log_2(n_1n_2\cdots n_t)\,.$$ 
\end{lem}
\begin{proof}
By induction, we can assume that $n_i = \ell_i + 1$ for exactly one choice of $i$, and $n_j = \ell_j$ for all $j \neq i$.  Note that we must have $n_i \geq 3$, and $\prod_{j \in \{1,\ldots,t\} \cut \{i\}} (n_j - 1)= \prod_{j \in \{1,\ldots,t\} \cut \{i\}} (\ell_j - 1) \geq 2$ by our hypotheses on $n_1,\ldots,n_t$. Thus, we have
\begin{align*}
    \prod_{k=1}^t (n_k - 1) - 2\log_2\left(\prod_{k=1}^t n_k\right) &= \left(\frac{n_i - 1}{\ell_i - 1}\right)\prod_{k=1}^t (\ell_k - 1) -\left( 2\log_2\left(\prod_{k=1}^t \ell_k\right) + 2\log_2\left(\frac{n_i}{\ell_i}\right)\right)\\
    &= \left(\frac{1}{\ell_i} + 1\right)\prod_{k=1}^t (\ell_k - 1) -  2\log_2\left(\prod_{k=1}^t \ell_k\right) - 2\log_2\left(\frac{\ell_i + 1}{\ell_i}\right)\\
    &\geq \prod_{j \in \{1,\ldots,t\} \cut \{i\}} (\ell_j - 1) - 2\log_2\left(\frac{\ell_i + 1}{\ell_i}\right)\\
    &\geq 2 - \log_2\left(\frac{\ell_i + 1}{\ell_i}\right)\\
    &> 0\,.
\end{align*}
    Observe that the last inequality follows since $\frac{\ell_i + 1}{\ell_i} < 2$, hence $\log_2\left(\frac{\ell_i + 1}{\ell_i}\right) < 1$.
\end{proof}

\begin{prop}\label{prop: direct exceptions}
Suppose we have positive integers $2 \leq n_1 \leq n_2 \leq \cdots \leq n_t$ with $n_2 \geq 3$ and $t \geq 2$; and additionally if $t = 2$, then $n_1 \geq 3$.  Then 
$$(n_1 - 1)(n_2 - 1)\cdots (n_t - 1) \geq 2\log_2(n_1n_2\cdots n_t)$$
unless $t = 2$ with
$$(n_1,n_2) \in \{(3,3), (3,4)\}\,,$$
or $t = 3$ with
$$(n_1,n_2,n_3) \in \{(2,3,3),(2,3,4),(2,3,5),(2,3,6),(2,4,4),(3,3,3)\}\,,$$
or $t = 4$ with 
$$(n_1,n_2,n_3,n_4) \in \{(2,3,3,3),(2,3,3,4)\}\,.$$
\end{prop}
\begin{proof}
This follows Lemma \ref{lem: cart product ineq} and basic case checking.
\end{proof}

Note that the inequalities appearing in Lemma \ref{lem: cart product ineq} and Proposition \ref{prop: direct exceptions} is precisely the minimum degree condition of Theorem \ref{thm:mindeg}.  We now combine these results to obtain the main theorem of this section.

\begin{thm}
The domination polynomials of all connected direct products of complete graphs are unimodal.
\end{thm}
\begin{proof}
Lemma \ref{lem:regular} implies that all graphs of the form $K_2 \times K_n$ have unimodal domination polynomials.  Combining Proposition \ref{prop: direct exceptions} and Theorem \ref{thm:mindeg}, we find that the domination polynomials of all connected direct products of complete graphs are unimodal except possibly the ten exceptional cases listed in Proposition \ref{prop: direct exceptions}.  It is straightforward to compute the domination polynomials in these ten remaining cases and see that the polynomials are indeed unimodal.
\end{proof}

\subsection{Cartesian Products of Complete Graphs}\label{sec:cartesian}

Motivated by our progress in understanding the direct product of complete graphs, we consider the domination polynomials of Cartesian products of complete graphs. We briefly show that unimodality holds for the Cartesian product of two complete graphs (known as a \textit{rook graph}), followed by a discussion of why higher order products are more difficult to handle than in the direct product case. 

\begin{thm}\label{thm:cartesian}
$D(K_m \cart K_n,x)$ is unimodal for all $m,n \in \N$.
\end{thm}
\begin{proof}
Without loss of generality we can assume $m \leq n$, since $K_m \cart K_n \cong K_n \cart K_m$.  When $m = 1$, $D(K_m \cart K_n,x)$ is unimodal as its coefficients are merely the binomial coefficients. When $m = 2$, we notice that $K_2 \cart K_n$ is an $n$-regular graph on $2n$ vertices. Therefore, by Lemma \ref{lem:regular}, $D(K_2 \cart K_n,x)$ is unimodal.  For $m \geq 3$, we can divide the set of remaining products of two complete graphs into three cases: 
\begin{enumerate}
    \item $m=3$
    \item $4\leq m\leq n \leq 7$
    \item $m \geq 4, n\geq 8$.
\end{enumerate}

Let us consider $m = 3$. We see for $n \geq 9$, $D(K_m \cart K_n,x)$ is unimodal by Theorem \ref{thm:mindeg}.  The cases for $n \leq 8$ can be checked explicitly. Similarly, the ten cases for $4 \leq m\leq n \leq 7$ can be checked explicitly. We see $D(K_6 \cart  K_7,x)$ and $D(K_7 \cart K_7,x)$ are unimodal by Theorem $\ref{thm:mindeg}$. The remaining eight cases can be verified by computer.

Therefore, it only remains to show that $D(K_m \cart  K_n,x)$ is unimodal for $m \geq 4, n\geq 8$. Recall Theorem \ref{thm:mindeg} proves that if $\delta_G \geq 2\log_2(|V|)$, then $D(G,x)$ is unimodal. So, we see for $m \geq 4$, $n \geq 8$, the minimum degree is given by
$$
    \delta_{K_m \cart K_n} = (m-1) + (n-1) \geq 2\log_2(m) + n-2 \geq 2\log_2(m) + 2\log_2(n) = 2\log_2(mn)\,.
$$
\end{proof}

We now discuss why higher order products are more infeasible than in the direct product case.  For a bounded number of complete graphs included in the Cartesian product, there is a finite and increasing number of cases that do not satisfy the hypotheses of Theorem \ref{thm:mindeg}.  Thus, while one could iterate this procedure for Cartesian products of three or more complete graphs, the number of exceptional cases, i.e., those graphs not covered by Theorem \ref{thm:mindeg}, quickly becomes computationally challenging to handle. Already for a Cartesian product of three complete graphs, there are $72$ exceptional cases, and $221$ for a product of four complete graphs.  Contrast this to the direct product case, where the minimum degree is $(n_1-1)(n_2-1)\cdots (n_t-1)$, which in general is much larger than the minimum degree $(n_1-1) + (n_2 - 1) + \cdots + (n_t - 1)$ in the Cartesian setting.  As we saw in Proposition \ref{prop: direct exceptions}, the hypotheses of Theorem \ref{thm:mindeg} were satisfied for all connected direct products of five or more complete graphs.  

One particular special case that may be interesting to investigate is the domination polynomial of the hypercube graph, defined recursively by $Q_n := K_2 \cart Q_{n-1}$ with $Q_1 := K_2$.  The problem of determining the domination number of $Q_n$ is a fundamental problem in coding theory and is only known in special cases, namely for $n \leq 9$ and $n = 2^k - 1$ or $2^k$, see \cite{HL, AHK}.  Thus, precisely determining the domination polynomials of hypercubes is very difficult, but it may still be possible to further understand the behavior of their coefficients.  

\section{Unimodality and Upper Domination}\label{sec:upper dom}

In this section, we use a Hall-type argument to show that for sufficiently large $i$, depending on the upper domination number, the coefficients $d_i(G)$ are non-increasing.

A subset $M \subseteq E(G)$ is a \emph{matching} if no two edges in $M$ share a vertex.  A matching $M$ \emph{saturates} a set of vertices $U \subseteq V(G)$ if every vertex of $U$ is incident to an edge in $M$.  We make use of Hall's classic theorem on saturated matchings in bipartite graphs.

\begin{thm}[Hall's Theorem, \protect{see, e.g., \cite[Theorem 5.2]{BM}}]\label{thm: Hall}
Let $G$ be a bipartite graph with bipartition $(X,Y)$.  Then $G$ contains a matching that saturates every vertex in $X$ if and only if, for all $U \subseteq X$, we have $|N(U)| \geq |U|$.
\end{thm}

Similar to the construction used by Alikhani and Peng to prove Theorem \ref{thm: first half}, we define a bipartite graph with vertices labeled by dominating sets of $G$.  This construction admits a saturated matching provided the dominating sets are sufficiently large, with the lower bound depending on the upper domination number. Recall that the \textit{upper domination number} $\Gamma(G)$ of a graph $G$ is the maximum size of a minimal dominating set. 

\begin{thm}\label{thm: upper dom}
Let $G$ be a graph on $n$ vertices.  If $i \geq \frac{n+\Gamma(G) - 1}{2}$, then $d_i(G) \geq d_{i+1}(G)$.  
\end{thm}
\begin{proof}
Let $A$ be a set of vertices labeled by the dominating sets of size $i+1$ in $G$, and let $B$ be a set of vertices labeled by the dominating sets of size $i$.  We construct a bipartite graph on $A \cup B$ by connecting two vertices if the corresponding dominating set of size $i$ is contained in the corresponding dominating set of size $i+1$.  Note that every vertex of $B$ has degree $n - i$.  Moreover, since a dominating set of size $i$ in $G$ contains a minimal dominating set of size at most $\Gamma(G)$, each vertex in $A$ has degree at least $i + 1 - \Gamma(G)$.  

We want to show that $A$ admits a saturated perfect matching into $B$.  By Theorem \ref{thm: Hall}, it is sufficient to show that $|N(T)| \geq |T|$ for every $T \subseteq A$.  Suppose that there exists $U \subseteq A$ such that $|N(U)| < |U|$.  Since each vertex of $U$ has degree at least $i + 1 - \Gamma(G)$, there must exist a vertex of $N(U)$ with degree greater than $i + 1 - \Gamma(G)$ in the induced subgraph on $U \cup N[U]$.  That is, we must have $i + 1 - \Gamma(G) < n - i$.  By our choice of $i$, we reach a contradiction.  Thus, $A$ admits a saturated perfect matching, so $|A| \leq |B|$, i.e., $d_{i+1}(G) \leq d_i(G)$.
\end{proof}

In particular, we obtain the following corollary.

\begin{cor}\label{cor: low upper dom}
Let $G$ be a graph on $n$ vertices.  If $n$ is odd and $\Gamma(G) \leq 4$, or if $n$ is even and $\Gamma(G) \leq 3$, then $D(G,x)$ is unimodal with a mode at $\left\lceil{\frac{n}{2}}\right\rceil$ or $\left\lceil{\frac{n}{2}}\right\rceil + 1$.
\end{cor}
\begin{proof}
Suppose first that $n$ is odd.  From Theorem \ref{thm: first half}, we have that 
$$d_1(G) \leq d_2(G) \leq \cdots \leq d_{\frac{n+1}{2}}(G)\,.$$
By Theorem \ref{thm: upper dom}, we furthermore have
$$d_{\frac{n+3}{2}}(G) \geq d_{\frac{n+5}{2}}(G) \geq \cdots \geq d_n(G)\,.$$
Thus, regardless of the relation between $d_{\frac{n+1}{2}}(G)$ and $d_{\frac{n+3}{2}}(G)$, the resulting sequence is unimodal.  The argument for $n$ even is analogous, but instead the only relation we cannot determine is between $d_{\frac{n}{2}}(G)$ and $d_{\frac{n+2}{2}}(G)$.
\end{proof}

\begin{rem}
Note that we could use Corollary \ref{cor: low upper dom} to show that the lollipop graphs $L_{m,n}$ have unimodal domination polynomials for $n \leq 3$.  However, we need to know additional information about the location of modes of these graphs in order to induct on $n$ with Theorem \ref{thm: BeBr poly recursion}, hence we handled this case more explicitly in Section \ref{sec:lollipop}.
\end{rem}

\section{Graphs with Universal Vertices}\label{sec:universal}
In this section, we investigate properties of the domination polynomial of a graph with many universal vertices. Let $G$ be a graph on $n$ vertices. A \emph{universal vertex} is a vertex whose closed neighborhood is the entire graph, i.e., a vertex of degree $n-1$. Since any set of vertices containing a universal vertex is necessarily dominating, we can place lower bounds on the coefficients $d_i(G)$ in terms of the number of universal vertices. This allows us to prove unimodality for graphs with sufficiently many universal vertices, in particular when the number of universal vertices is at least $\log_2(n) -1$.

We use the following result of Beaton and Brown to show that the coefficients of the domination polynomial for a graph with one or more universal vertices are non-increasing past a certain threshold. For a graph on $n$ vertices, let $r_k(G)$ be the proportion of subsets of $V(G)$ of size $k$ that are dominating. That is, $r_k(G) = \frac{d_k(G)}{{n \choose k}}.$

\begin{lem}\emph{(\cite[Lemma 3.1]{BeBr})}\label{lem: domination ratio}
Let $G$ be a graph on $n$ vertices, and $k\geq \frac{n}{2}$.  If $r_k(G)\geq \frac{n-k}{k+1}$ then $d_{i+1}(G)\leq d_i(G)$ for all $i\geq k$. In particular, if $k=\left\lceil\frac{n}{2}\right\rceil$ then $D(G,x)$ is unimodal with a mode at $\left\lceil\frac{n}{2}\right\rceil$.
\end{lem}

\begin{thm}\label{thm: universal decreasing}
Let $G$ be a graph on $n$ vertices, $m$ of which are universal.  Let $w_m$ be the smallest positive real root of $f_m(x) = x^{m+1}-x^m - 2x + 1$.  Then $d_i(G) \geq d_{i+1}(G)$ for all $i \geq (1-w_m)n$.
\end{thm}
\begin{proof}
Suppose that we have $d_{i+1}(G) > d_{i}(G)$ for some $i \geq \frac{n}{2}$, and we wish to show that $i < (1-w_m)n$.  Since any non-dominating set of size $i$ must not contain any of the $m$ universal vertices, we have
$$d_i(G) \geq {n \choose i} - {n - m \choose i}\,.$$
By Lemma \ref{lem: domination ratio}, we must have
$$\frac{d_{i}(G)}{{n \choose i}} < \frac{n-i}{i+1}\,.$$
Combining the above two inequalities yields
$${n \choose i} - {n - m \choose i} \leq d_i(G) \leq \frac{n-i}{i+1}{n \choose i}\,.$$
We can then divide both sides by ${n \choose i}$ and rearrange to obtain
$$\frac{n-i}{i+1} + \prod_{j=0}^{m-1} \frac{n-i -j}{n-j} > 1\,.$$
Since $\frac{n-i-j}{n-j} < \frac{n-i}{n}$ for any $j > 0$, we furthermore have
$$\frac{n-i}{i} + \left(\frac{n-i}{n}\right)^m \geq \frac{n-i}{i+1} + \left(\frac{n-i}{n}\right)^m > 1\,.$$
Setting $i = (1-\kappa)n$, we see that this condition is equivalent to 
$$\frac{\kappa n}{(1-\kappa) n} + \left (\frac{\kappa n}{n}\right)^m - 1 = \frac{\kappa}{1-\kappa} + \kappa ^m  - 1 > 0\,.$$
Multiplying by $\kappa - 1$ (which is necessarily negative) yields
$$\kappa^{m+1}-\kappa^m - 2\kappa + 1 < 0.$$
Since $f(0) > 0$ and $f(\kappa) < 0,$ by the above inequality, we must have that $\kappa > w_m$.  Hence, we can conclude $i = (1-\kappa)n < (1-w_m)n$, as desired.
\end{proof}

We can now use a bound on the smallest positive root of $f_m(x)$, the proof of which is included at the end of this section, to show that that certain domination polynomials are unimodal.

\begin{lem}\label{lem: root estimate}
For $m \in \N$, the smallest positive real root of $f_m(x) = x^{m+1} - x^m -2x + 1$ lies in the range $\left(\frac{1}{2} - \frac{1}{2^{m+1}}, \frac{1}{2} \right)$.
\end{lem}

When the above bounds $w_m$ are sufficiently close to $\frac{1}{2}$ with respect to the number of vertices $n$, we can conclude that the latter half of the domination coefficients is non-increasing.

\begin{thm}\label{thm: unimodal many universal}
If $G$ is a graph on $n$ vertices having at least $\log_2(n) - 1$ universal vertices, then $D(G,x)$ is unimodal with a mode at $\left\lceil\frac{n}{2}\right\rceil$ or $\left\lceil\frac{n}{2}\right\rceil + 1$.
\end{thm}
\begin{proof}
Let $m$ denote the number of universal vertices of $G$.  By Theorem \ref{thm: universal decreasing}, we have 
\\$d_{i}(G) \geq d_{i+1}(G)$ for $i \geq (1-w_m)n$.  Using the estimate $\frac{1}{2} - w_m < \frac{1}{2^{m+1}}$ from Lemma \ref{lem: root estimate} along with our assumption that $m \geq \log_2(n)-1$, this implies 
$$(1-w_m)n < \left(\frac{1}{2} + \frac{1}{2^{m+1}}\right)n = \frac{n}{2} + \frac{n}{2^{m+1}} \leq \frac{n}{2} + 1\,.$$
Therefore, we have $d_i(G) \geq d_{i+1}(G)$ for $i \geq \frac{n}{2} + 1$.  Since, by Theorem \ref{thm: first half}, we also have $d_i(G) \leq d_{i+1}(G)$ for $i < \frac{n}{2}$, we can conclude that $d_i(G)$ is unimodal with a mode as claimed.
\end{proof}

\begin{rem}
Note that the minimum degree of a graph satisfying the hypotheses of Theorem \ref{thm: unimodal many universal} can be as low as $\log_2(n) - 1$.  Thus, these graphs may have minimum degree as low as half of the lower threshold required to apply Theorem \ref{thm:mindeg} by Beaton and Brown, which assumes that the minimum degree is at least $2\log_2(n)$.  
\end{rem}

Lastly, we prove the claimed bounds on $w_m$, the smallest positive root of the polynomial $f_m(x) = x^{m+1}-x^m - 2x + 1$.

\begin{proof}[Proof of Lemma \ref{lem: root estimate}]
By Descartes' Rule of Signs, this polynomial can have at most two positive real roots.  Moreover, since $f_m(0), f_m(2) > 0$ and $f_m\left(\frac{1}{2}\right) < 0$, $f_m$ must have one root in the interval $\left(0,\frac{1}{2}\right)$, which we denote by $w_m$, and another in the interval $\left(\frac{1}{2},2\right)$. 

We expect $w_m$ to approach $\frac{1}{2}$ from below as $m$ increases.  So we let $\varepsilon_m = \frac{1}{2} - w_m > 0$ and estimate $\varepsilon_m$.  Using the fact that $(1+ \frac{1}{y})^y > e$ for all $y < -1$, we have
$$(w_m)^m = \left(\frac{1-2\varepsilon_m}{2}\right)^m = \frac{(1-2\varepsilon_m)^{-\frac{1}{2\varepsilon_m}(-2\varepsilon_mm)}}{2^m} \leq \frac{e^{-2\varepsilon_mm}}{2^m}\,.$$
Now, additionally using the hypothesis that $2w_m - 1 = (w_m)^m(w_m-1)$, we have
$$2\left(\frac{1}{2}-\varepsilon_m\right) - 1 \geq \frac{e^{-2\varepsilon_m}}{2^m}\left(\left(\frac{1}{2}-\varepsilon_m\right)-1\right)\,.$$
In particular, recalling that $\varepsilon_m < \frac{1}{2}$ we obtain
$$\varepsilon_m \leq \frac{e^{-2\varepsilon_mm}}{2^{m+1}}\left(\frac{1}{2} + \varepsilon_m\right) \leq \frac{1}{2^{m+1}}\left(\frac{1}{2} + \varepsilon_m\right) = \frac{1}{2^{m+2}} + \frac{\varepsilon_m}{2^{m+1}} < \frac{1}{2^{m+1}}\,.$$
This yields the corresponding lower bound for $w_m$.
\end{proof}

\section{Further Directions}\label{sec: further directions}
As discussed in Section \ref{sec:spiders}, the unimodality of the domination polynomials for families of graphs with low minimum degree is of particular interest.  We show that spiders with at most $400$ legs have unimodal domination polynomials, so a natural extension would be showing this result for an arbitrary number of legs.  If Conjecture \ref{conj: spider modes} holds, this would be sufficient to prove unimodality for all spider graphs. Part of this conjecture asserts that the operation of leaf deletion yields a domination polynomial having a mode of distance at most $1$ from a mode of the original domination polynomial.  In fact, we have not yet found a counterexample to this phenomenon for trees, leading to the following question. 

\begin{ques}
    Does there exist a tree $T$ and a leaf $v \in V(T)$ such that $D(T)$ and $D(T \cut v)$ are both unimodal but have no modes of distance at most $1$?
\end{ques}

In Subsection \ref{sec:cartesian}, we show that the Cartesian product of any two complete graphs is unimodal with a mode at $\left\lceil\frac{n}{2}\right\rceil$ using Theorem \ref{thm:mindeg}.  If one considers Cartesian products of arbitrarily many complete graphs, it is straightforward to see that there are infinitely many graphs that do not satisfy the hypotheses of Theorem \ref{thm:mindeg}, unlike in the direct product case.  It would be interesting to develop new techniques for showing unimodality for this family, which may generalize to higher order products.

We examined graphs with many vertices of the highest possible degree in Section \ref{sec:universal}.  It seems likely that the results in this section could be generalized for graphs with many vertices of the near highest degree, i.e., for graphs with $m$ vertices of degree $n - k$ for fixed $k$ and sufficiently large $n$.  In this setting, it would be interesting to determine the threshold for $m$ (in terms of $n$ and $k$) after which one could determine that the domination polynomial is unimodal with a mode at $\left\lceil \frac{n}{2}\right\rceil$.

\section{Acknowledgements}
The authors thank Joe Gallian for his tireless efforts leading the Duluth REU community.   They furthermore thank Colin Defant, Milan Haiman, Noah Kravitz, and Yelena Mandelshtam for helpful suggestions and insightful questions.  This research was conducted at the University of Minnesota Duluth Mathematics REU and was supported, in part, by NSF-DMS Grant 1949884 and NSA Grant H98230-20-1-0009.  The first author was also supported by a Marshall scholarship, and additional support for the second author was provided by the Mitya.

\end{document}